\documentclass{amsart}
\title[Interpolation of hypergeometric ratios in a global field]{Interpolation of hypergeometric ratios in a global field of positive characteristic}
\author{Greg W.\ Anderson}
\date{June 26, 2006}
\email{gwanders@math.edu}
\address{School of Mathematics, University of Minnesota, Minneapolis, MN 55455}
\newcommand{\Abel}{{\mathrm{Abel}}}
\newcommand{\one}{{\mathbf{1}}}
\DeclareMathOperator{\Aut}{{\mathrm{Aut}}}
\newcommand{\UUU}{{\mathcal{U}}}
\DeclareMathOperator{\supp}{{\mathrm{supp}}}
\newcommand{\Mbold}{{\mathbf{M}}}
\newcommand{\Dbold}{{\mathbf{D}}}

\newcommand{\Kbold}{{\mathbf{K}}}
\DeclareMathOperator{\ord}{{\mathrm{ord}}}
\DeclareMathOperator{\Div}{{\mathrm{Div}}}
\DeclareMathOperator{\Gal}{{\mathrm{Gal}}}
\DeclareMathOperator{\sep}{{\mathrm{sep}}}
\DeclareMathOperator{\perf}{{\mathrm{perf}}}
\DeclareMathOperator{\ab}{{\mathrm{ab}}}
\DeclareMathOperator{\adeles}{{\mathbb A}}

\DeclareMathOperator{\trace}{{\mathrm{trace}}}

\DeclareMathOperator{\Moore}{{\mathrm{Moore}}}

\DeclareMathOperator{\PP}{{\mathbb{P}}}

\newcommand{\iso}{\xrightarrow{\sim}}
\newcommand{\OO}{{\mathcal{O}}}
\DeclareMathOperator{\Spec}{{\mathrm{Spec}}}

\DeclareMathOperator{\Res}{{\mathrm{Res}}}
\DeclareMathOperator{\RES}{{\mathrm{RES}}}
\newcommand{\FF}{{\mathbb{F}}}
\newcommand{\ZZ}{{\mathbb{Z}}}

\newcommand{\norm}[1]{{\Vert #1 \Vert}}
\DeclareMathOperator{\Frob}{\mathrm{Frob}}
\DeclareMathOperator{\Hyp}{{\mathrm{Hyp}}}
\newtheorem{Theorem}[subsection]{Theorem}
\newtheorem{Lemma}[subsubsection]{Lemma}
\newtheorem{Proposition}[subsection]{Proposition}
\newcommand{\AUT}[1]{{\Aut(\overline{\FF_q(#1)}/
\overline{\FF}_q)}}

\begin{document}
\begin{abstract}
In connection with each global field of positive characteristic we exhibit many examples of two-variable algebraic functions possessing properties consistent with a conjectural refinement of the Stark conjecture in the function field case recently proposed by the author ({\tt math.NT/0407535}). Most notably, all examples are Coleman units.
 We obtain our results by studying rank one shtukas in which both zero and pole are generic, i.~e., shtukas not associated to any Drinfeld module. 
 \end{abstract}
\maketitle
\tableofcontents

\section{Introduction}
Our main result (Theorem~\ref{Theorem:MainResult} below) provides in connection with each global field of positive characteristic  many examples of two-variable algebraic functions with at least some  properties predicted by the author's conjecture \cite[Conj.\ 9.5]{AndersonStirling}. 
Most notably, each example is a Coleman unit.
Furthermore, each example figures in an  interpolation formula  in which the hypergeometric ratios mentioned in the title of the paper appear on the right side. 
The notion of Coleman unit, which was inspired by Coleman's remarkable paper \cite{Coleman}, was introduced in \cite{AndersonStirling}
and is reviewed in \S\ref{subsection:ColemanUnits} below.
The notion of interpolation formula can be traced back to papers of Thakur, especially \cite{ThakurAnnals} and \cite{ThakurIMRN}; roughly speaking, in such a formula a Frobenius endomorphism appears on the left side raised to a variable power.
The notion of hypergeometric ratio, which is  a specialization of the notion of Catalan symbol introduced in \cite{AndersonStirling}, is defined in \S\ref{subsection:HypergeometricRatios} below.

Our constructions are based on the study of rank one shtukas in a relatively elementary setting similar to that of Thakur's paper \cite{Thakur}.  The Coleman units we produce come into existence as  invariants naturally attached to shtukas.
But the new twist here in comparison to \cite{Thakur}
is that our shtukas have both generic zero and generic pole, and hence are not attached to any Drinfeld module.

Ultimately an analysis of the examples constructed here with tools developed in \cite{AndersonStirling} and \cite{AndersonLocalStirling} 
 yields a proof of \cite[Conj.\ 9.5]{AndersonStirling},
but because the bookkeeping needed to complete that proof is heavy and lengthy, we will provide the details on another occasion. Here we will just focus on the construction of Coleman units satisfying  interpolation formulas.
The main point we want to make is that
the Coleman unit property follows naturally from a variant (Lemma~\ref{Lemma:ChiZero} below) of Drinfeld's powerful ``$\chi=0\Rightarrow h^0=h^1=0$'' lemma
\cite{Drinfeld} (see also \cite[p.\ 146]{Mumford}).

We consider this paper to be third in a series starting with \cite{AndersonStirling} and \cite{AndersonLocalStirling}, and accordingly we recommend that the reader scan the introductions of those papers for background, motivation, and further references. (The introduction to \cite{ABP} might also be helpful.) 
But no detailed familiarity with
\cite{AndersonStirling} and \cite{AndersonLocalStirling}
is assumed here. This paper is largely
independent of its predescessors.

\section{Formulation of the main result}\label{section:Formulation}
\subsection{Basic setting and notation}

\subsubsection{The curve $X/\FF_q$}
Let $X/\FF_q$ be a smooth projective geometrically connected curve of genus $g$, where the base $\FF_q$ is a field of $q<\infty$ elements. The curve $X/\FF_q$
remains fixed throughout the paper. We denote the function field of $X$ by $\FF_q(X)$.
We use standard notation for coherent sheaves
and cohomology on $X$.  We will deal with no sheaves more complicated than invertible sheaves and  their quotients.

\subsubsection{Moore determinants}\label{subsubsection:Moore}
Put
$$\Moore(x_1,\dots,x_n)=
\left|\begin{array}{lll}
x_1^{q^{n-1}}&\dots&x_n^{q^{n-1}}\\
\vdots&&\vdots\\
x_1^{q^0}&\dots&x_n^{q^0}
\end{array}\right|\in \FF_q[x_1,\dots,x_n]$$
where $x_1,\dots,x_n$ are independent variables.
Recall the {\em Moore determinant identity}:
$$\Moore(x_1,\dots,x_n)
=
\prod_{k=1}^n
\left(\prod_{a_{k+1}\in \FF_q}\cdots
\prod_{a_{n}\in \FF_q}
\left(x_k+a_{k+1}x_{k+1}+\cdots+a_nx_n\right)\right).$$
\subsubsection{Residues}\label{subsubsection:RES}
Given an effective divisor $D$ of $X$
and a meromorphic differential $\omega$ on $X$, we define
$\RES_D \omega$ to be the sum of terms $\trace_{\FF_x/\FF_q}\Res_x\omega$
extended over closed points $x$ of $X$ in the support of $D$, where $\FF_x$ is the residue field at $x$ and $\Res_x\omega\in \FF_x$ is the residue of $\omega$ at $x$. Note that $\RES_D$ induces a perfect $\FF_q$-bilinear pairing $H^0(\OO_X(D)/\OO_X)\times 
H^0(\Omega_{X/\FF_q}/\Omega_{X/\FF_q}(-D))\rightarrow \FF_q$.

\subsubsection{Generalized divisor classes}\label{subsubsection:Gen1}
Given an effective divisor $D$ of $X$
and a  nonzero meromorphic function $f$ on $X$, we write $f\vert_D\equiv 1$
if $f$ is regular in a neighborhood of $D$ 
and its restriction $f\vert_D$ to the closed subscheme $D$ is identically equal to $1$, in which case we also
say that the divisor $(f)$ is {\em principal to the conductor $D$}.
Given an effective divisor $D$ of $X$
and divisors $E_1$ and $E_2$ of $X$ supported away from $D$, we say that $E_1$ and $E_2$ belong to the same {\em generalized divisor class} of {\em conductor $D$} 
and we write $E_1\sim_DE_2$ if $E_1-E_2$ is principal to the conductor $D$.

\subsubsection{Miscellaneous}
Let $A^\times$ denote the multiplicative group of a ring $A$ with unit. 

\subsection{Apparatus from class field theory}

\subsubsection{The id\`{e}le group of $X$}
Let $\adeles_X$ (resp., $\adeles_X^\times$) be the ad\`{e}le ring 
(resp., id\`{e}le group) of $X$.
We identify $\FF_q(X)^\times$
with the diagonal subgroup of $\adeles_X^\times$, as usual.
Let 
$$\norm{\cdot}:\adeles_X^\times \rightarrow q^\ZZ$$
be the id\`{e}le norm homomorphism.
To each id\`{e}le $a\in \adeles_X^\times$ we associate a divisor 
$$\Div a=\sum_x (\ord_x a)x,$$
where the sum is extended over closed points $x$ of $X$, and $\ord_x a$
denotes the order of vanishing of $a$ at $x$. The rule $\Div$ extends the usual rule for associating a divisor to a meromorphic function on $X$.
Note that 
$$-\deg  \Div a=\log_q\norm{a}$$ for all $a\in \adeles_X^\times$. Given an effective divisor $D$ of $X$ and $a\in \adeles_X^\times$,
we say that $a$ is {\em supported away from $D$} if for every
closed point $x$ in the support of $D$ we have 
$$\ord_x (a-1)\geq \ord_x D,$$
in which case the divisor $\Div a$ is 
also supported away from $D$.

\subsubsection{The reciprocity law homomorphism}
Let $\overline{\FF_q(X)}$ be an algebraic closure of $\FF_q(X)$.
Let $\FF_q(X)^{\ab}$ be the abelian closure of $\FF_q(X)$ in $\overline{\FF_q(X)}$.
Let $\FF_q(X)_{\perf}$ (resp., $\FF_q(X)^{\ab}_{\perf}$) be the closure of
$\FF_q(X)$ (resp., $\FF_q(X)^{\ab}$) in $\overline{\FF_q(X)}$
under the extraction of $q^{th}$ roots. 
We define
$$\rho:\adeles_X^\times\rightarrow \Gal(\FF_q(X)^{\ab}/\FF_q(X))=
\Gal(\FF_q(X)^{\ab}_{\perf}/\FF_q(X)_{\perf})$$
to be the reciprocity law homomorphism of global class field theory,
``renormalized'' in the fashion of \cite{TateBackground} so that
$$\rho(a) C=C^{\norm{a}}$$ holds for every $C$ belonging
to the algebraic closure $\overline{\FF}_q$ of $\FF_q$ in 
$\overline{\FF_q(X)}$ and $a\in \adeles_X^\times$.  We define
$$\rho^*:\adeles_X^\times\rightarrow
\Aut(\FF_q(X)_{\perf}^{\ab}/\overline{\FF}_q)$$
by the rule
$$\rho^*(a)x=(\rho(a)^{-1}x)^{\norm{a}}$$
for all $x\in \FF_q(X)_{\perf}^{\ab}$ and $a\in \adeles_X^\times$. 
The homomorphism $\rho^*$ actually plays a more important role in this paper than does $\rho$.

\subsubsection{The homomorphism $r_D$}
\label{subsubsection:IdeleQuotient}
Let $D$ be an effective divisor of $X$.
Let $\UUU_D\subset \adeles_X^\times$ be the open compact subgroup consisting
of id\`{e}les $a$ such that for all closed points $x\in X$,
if $x$ is (resp., is not) in the support of $D$, then
$\ord_x(a-1)\geq \ord_x D$ (resp., $\ord_x D=0$).
There is a unique exact sequence
\begin{equation}\label{equation:rD}
1\rightarrow \FF_q(X)^\times\UUU_D\subset \adeles_X^\times
\xrightarrow{r_D} \left(\begin{array}{l}
\mbox{generalized divisor class}\\
\mbox{group of conductor $D$}
\end{array}\right)\rightarrow 0
\end{equation}
such that 
$$r_D(a)=\left(\begin{array}{l}
\mbox{generalized divisor class }\\
\mbox{of $-\Div a$ of conductor $D$}
\end{array}\right)
$$
for every id\`{e}le $a\in \adeles_X^\times$ supported away from $D$.

\subsubsection{Remark}\label{subsubsection:Cancellation}
Let $D$ be an effective divisor of $X$.
Let $K/\FF_q(X)$ be a finite abelian extension of conductor dividing $D$.
Let $x$ be a closed point of $X$ not in the support of $D$
and hence unramified in $K/\FF_q(X)$.
Let $\sigma_x\in \Gal(K/\FF_q(X))$ be the arithmetic Frobenius element at $x$, i.~e., the traditional value of the Artin symbol $(x,K/\FF_q(X))$.
 Let $a\in \adeles_X^\times$ be such that $r_D(a)=x$.
Then we have $\rho(a)\vert_K=\sigma_x$. 
In a nutshell: the minus sign intervening in the definition of $r_D$
cancels the renormalization of $\rho$.

\subsection{Hypergeometric ratios}\label{subsection:HypergeometricRatios}
We introduce  a notion which is actually a specialization
of the notion of Catalan symbol introduced in \cite{AndersonStirling}.

 \subsubsection{Definition (high degree case)}
 Let $D$ be a nonzero effective divisor of $X$.
Let $E$ be a divisor of $X$ supported away from $D$.
Assume that $\deg E>2g-2$, in which case $H^1(\OO_X(E))=0$, 
and hence the sequence
$$0\rightarrow H^0(\OO_X(E))\rightarrow H^0(\OO_X(E+D))\rightarrow
H^0(\OO_X(D)/\OO_X)\rightarrow 0$$
is exact. Let nonzero $\alpha,\beta\in H^0(\OO_X(D)/\OO_X)$ be given, along with liftings $\tilde{\alpha},\tilde{\beta}\in H^0(\OO_X(E+D))$, respectively,
via the exact sequence above.  In this situation we define
$$\Hyp_D(\alpha,\beta,E)= 
\prod_{
e\in H^0(\OO_X(E))}\frac{\tilde{\alpha}+e}
{\tilde{\beta}+e}\in \FF_q(X)^\times,
$$
which is independent of the choice of liftings $\tilde{\alpha}$ and $\tilde{\beta}$.   We call $\Hyp_D(\alpha,\beta,E)$ a {\em hypergeometric ratio}.  Note that $\Hyp_D(\alpha,\beta,E)$ depends only on the generalized divisor class of $E$ of conductor $D$. More generally,
we have
\begin{equation}\label{equation:HypInv}
\Hyp_D(\alpha,\beta,E+(f))=\Hyp_D((f\vert_D)\alpha,(f\vert_D)\beta,E)
\end{equation}
for all $f\in \FF_q(X)^\times$ such that $(f)$ is supported away from $D$. We have
\begin{equation}\label{equation:MooreHyper}
\Hyp_D(\alpha,\beta,E)=\frac{\Moore(\tilde{\alpha},e_1,\dots,e_n)}
{\Moore(\tilde{\beta},e_1,\dots,e_n)}
\end{equation}
for every $\FF_q$-basis 
$$e_1,\dots,e_n\in H^0(\OO_X(E))\;\;\;(n=h^0(\OO_X(E))=\deg E-g+1),$$
whence follow the relations
\begin{equation}\label{equation:MooreHyperBis}
\Hyp_D(c\alpha,\beta,E)=\Hyp_D(\alpha,c^{-1}\beta,E)=c\Hyp_D(\alpha,\beta,E)
\end{equation}
for all $c\in \FF_q^\times$ and
\begin{equation}\label{equation:MooreHyperTer}
\Hyp_D(\alpha,\beta,E)=\Hyp_D(\alpha_1,\beta,E)+\Hyp_D(\alpha_2,\beta,E)
\end{equation}
for all decompositions $\alpha=\alpha_1+\alpha_2$ where
$\alpha_1,\alpha_2\in H^0(\OO_X(D)/\OO_X)$ are nonzero.

\subsubsection{Definition (low degree case)}
\label{subsubsection:LowDegree}
 As in the previous paragraph, let $D$ be a nonzero effective divisor of $X$
 and let $E$ be a divisor of $X$ supported away from $D$.
But this time let us assume that
$\deg E<-\deg D$, in which case $h^1(\Omega_{X/\FF_q}(-E-D))=0$ and hence the sequence
$$0\rightarrow H^0(\Omega_{X/\FF_q}(-E-D))
\rightarrow H^0(\Omega_{X/\FF_q}(-E))\rightarrow 
H^0(\Omega_{X/\FF_q}/\Omega_{X/\FF_q}(-D))
\rightarrow 0$$
is exact. Let nonzero $\alpha,\beta\in H^0(\OO_X(D)/\OO_X)$ be given, along with liftings $\tilde{\alpha}$ and $\tilde{\beta}$ to meromorphic functions on $X$, respectively.
In this situation we define
$$ 
\Hyp_D(\alpha,\beta,E)=\prod_{\begin{subarray}{c}
\omega\in H^0(\Omega_{X/\FF_q}(-E))\\
\RES_D(\omega \tilde{\beta})=1
\end{subarray}}\omega \bigg /
\prod_{\begin{subarray}{c}
\omega\in H^0(\Omega_{X/\FF_q}(-E))\\
\RES_D(\omega \tilde{\alpha})=1
\end{subarray}}\omega\;\in \FF_q(X)^\times,
$$
which is independent of the choice of liftings $\tilde{\alpha}$ and $\tilde{\beta}$. The ratio does indeed define a meromorphic function on $X$ because there are exactly $q^{g-2-\deg E}$ factors in 
the numerator, and an equal number of factors in the denominator. 
 Note that in the low degree case, just as in the high degree case, $\Hyp_D(\alpha,\beta,E)$ depends only on the generalized divisor class of $E$ to the conductor $D$, and furthermore satisfies (\ref{equation:HypInv}).
Trivially, formula
(\ref{equation:MooreHyperBis})
continues to hold.
Perhaps surprisingly, formula (\ref{equation:MooreHyperTer}) also continues to hold in the low degree case---this will follow from our main result,
and is anyhow easy to verify directly
using tricks discussed in \cite[\S3]{AndersonStirling}.

\subsubsection{Remark}\label{subsubsection:CatalanComparison}
This remark will not be needed to follow the main line of inquiry.
But it will be needed to make sense of later remarks.
Given a divisor $E$ of $X$, let us associate
to it an open compact  subgroup $[E]\subset\adeles_X$
by the rule
$$[E]=\{a\in \adeles_X\vert
\ord_x a+\ord_x E\geq 0\;\mbox{for all closed points $x\in X$}\}.$$
This rule has the property that
$$[E]\cap \FF_q(X)=H^0(X,\OO_X(E)).$$
Now fix a nonzero effective divisor $D$ of $X$
and 
$$\alpha,\beta\in H^0(\OO_X(D)/\OO_X)=[D]/[0].$$
Fix liftings
$$\tilde{\alpha},\tilde{\beta}\in [D]\subset \adeles_X,$$
respectively. Fix also a divisor $A_0$ of $X$
supported away from $D$ of degree $g-2$.
Given any subset $S\subset \adeles_X$, let $\one_S$ be the 
$\{0,1\}$-valued function on $\adeles_X$ taking the value $1$ on $S$ and $0$ elsewhere. It can be shown that
\begin{equation}\label{equation:CatalanSymbolSpecialization}
\left(\begin{array}{c}
a\\
\one_{\tilde{\alpha}+[A_0]}-
\one_{\tilde{\beta}+[A_0]}
\end{array}\right)=\Hyp_D(\alpha,\beta,A_0+r_D(a))^{\min(\norm{a},1)}
\end{equation}
for all $a\in \adeles_X^\times$ such that the right side is defined, where the object $(\begin{subarray}{c}\cdot\\\cdot\end{subarray})$ on the left is the {\em Catalan symbol}
defined in \cite{AndersonStirling}.
We omit the details of the comparison since we
wish to avoid introducing a lot of machinery of harmonic analysis which otherwise we will not be using.

\subsection{The ring $\Dbold$}
\subsubsection{Definitions}
Consider the ring
$$\Dbold=\overline{\FF_q(X)}\otimes_{\overline{\FF}_q}
\overline{\FF_q(X)}.$$
We define the {\em diagonal evaluation} homomorphism
$$(\varphi\mapsto \varphi\vert_\Delta):\Dbold\rightarrow
\overline{\FF_q(X)}$$
by the rule
$$(x\otimes y)\vert_\Delta=xy,$$
and correspondingly we define
$$\Delta=\ker\left(\varphi\mapsto \varphi\vert_\Delta\right)\subset \Dbold,$$
which is a maximal ideal of $\Dbold$.
For all $\theta_1,\theta_2\in \Aut(\overline{\FF_q(X)}/\FF_q)$ such that $$\theta_1\vert_{\overline{\FF}_q}=\theta_2\vert_{\overline{\FF}_q}$$
we define
$$\theta_1\otimes \theta_2:\Dbold\rightarrow\Dbold$$
by the rule
$$(\theta_1\otimes \theta_2)(x\otimes y)=(\theta_1x)\otimes (\theta_2 y).$$
In the case $(\theta_1,\theta_2)=(\mbox{identity automorphism},\theta)$ we write $\theta_1\otimes \theta_2=1\otimes \theta$.
\begin{Lemma}\label{Lemma:DGeography}
(i) The ring $\Dbold$ is a domain. (ii) Every nonzero ideal of $\Dbold$ is maximal.
(iii) The local ring $\Dbold_\Delta$ of $\Spec(\Dbold)$ at $\Delta$ is a nondiscrete valuation ring of rank one.
(iv) Every maximal ideal $\Mbold\subset \Dbold$ is of the form
$$\Mbold=\ker\left((\varphi\mapsto ((1\otimes \theta)\varphi)\vert_\Delta):
\Dbold\rightarrow \overline{\FF_q(X)}\right)$$
for unique $\theta\in \AUT{X}$.  
\end{Lemma}
\proof Let $L/\FF_q(X)$ be a finite subextension of $\overline{\FF_q(X)}/\FF_q(X)$ and put $\FF_\ell=L\cap\overline{\FF}_q$.
 Realize $L$ as the function field
$\FF_\ell(Y)$ of a smooth projective geometrically connected curve $Y/\FF_\ell$.  Given also a finite nonempty set $S$ of closed points of $Y$, let $$\Dbold_{L,S}=\overline{\FF_q(X)}\otimes_{\FF_\ell} H^0({Y\setminus S},\OO_Y).$$ 
The ring $\Dbold_{L,S}$ is the coordinate ring of an irreducible smooth affine curve
defined over the field $\overline{\FF_q(X)}$ and in particular is a Dedekind domain. Moreover, by the Nullstellensatz, the maximal ideals of $\Dbold_{L,S}$ correspond bijectively to $(\overline{\FF_q(X)}\otimes 1)$-linear homomorphisms $\Dbold_{L,S}\rightarrow
\overline{\FF_q(X)}$. Let $\Dbold_L$ be the limit over $S$ of $\Dbold_{L,S}$. Again $\Dbold_L$ is a Dedekind domain and maximal ideals of $\Dbold_L$
correspond bijectively to $(\overline{\FF_q(X)}\otimes 1)$-linear homomorphisms
$\Dbold_L\rightarrow\overline{\FF_q(X)}$. 
Given a tower $L_2/L_1/\FF_q(X)$ contained in $\overline{\FF_q(X)}/\FF_q(X)$ with $L_2/\FF_q(X)$ finite, the ring extension $\Dbold_{L_2}/\Dbold_{L_1}$ is finite flat, and moreover \'{e}tale if $L_2/L_1$ is separable.
The ring $\Dbold$ is the union of rings of the form $\sqrt[q^n]{\Dbold_{L}}$ with $L/\FF_q(X)$ ranging over finite separable subextensions of $\overline{\FF_q(X)}/\FF_q(X)$ and $n$ ranging over positive integers. The result follows
by passage to the limit on $L$ and $n$.
\qed

\subsubsection{Extensions}
For all $\theta_1,\theta_2\in\Aut(\overline{\FF_q(X)}/\FF_q)$ with a common restriction to $\overline{\FF}_q$ we extend the automorphism $\theta_1\otimes \theta_2$ of $\Dbold$ to the fraction field of
$\Dbold$ in the unique possible way. 
We extend diagonal evaluation to a homomorphism
$$(\varphi\mapsto \varphi\vert_\Delta):\Dbold_\Delta\rightarrow
\overline{\FF_q(X)}$$
in the unique possible way, and for convenience we set $\varphi\vert_\Delta=\infty$
for every $\varphi$ in the fraction field of $\Dbold$ which does not belong to $\Dbold_\Delta$. 

\begin{Lemma}\label{Lemma:UniquenessPrinciple}
Let $\varphi$ be an element of the fraction field of $\Dbold$ such that
for infinitely many integers $n$ there exists  $\theta\in \AUT{X}$ with the following two properties:
$\theta\vert_{\FF_q(X)_{\perf}}=(x\mapsto x^{q^n})$
and $((1\otimes \theta)\varphi)\vert_\Delta=0$. Then $\varphi=0$.
\end{Lemma}
\proof Notation as in the proof of Lemma~\ref{Lemma:DGeography}, the function $\varphi$
belongs to the fraction field of some Dedekind domain of the form $\Dbold_{L}$.
By hypothesis $\varphi$ has positive valuation at infinitely many
distinct maximal ideals of $\Dbold_{L}$, and hence vanishes identically.
\qed

\subsubsection{Critical automorphisms and their exponents}
Given $\theta\in \AUT{X}$, we say that $\theta$ is {\em critical}
if there exists $a\in \adeles_X^\times$
such that 
$$\theta\vert_{\FF_q(X)^{\ab}_{\perf}}=\rho^*(a),$$
in which case $a$ is uniquely determined by $\theta$,
and will be called the {\em exponent} of $\theta$.

\begin{Lemma}\label{Lemma:Slick}
Fix $\theta\in \AUT{X}$. The following properties are equivalent:
\begin{itemize}
\item $\theta$ is critical.
\item $\theta\vert_{\FF_q(X)_{\perf}}=(x\mapsto x^{q^n})$
for some integer $n$.
\end{itemize}
\end{Lemma}
\proof The first property trivially implies the second.
The second property granted, the automorphism $\sqrt[q^n]{\theta}$
fixes every element of $\FF_q(X)_{\perf}$, stabilizes $\FF_q(X)^{\ab}_{\perf}$,
and restricts on $\overline{\FF}_q$ to an integer power of the $q^{th}$ power Frobenius automorphism. But then $\sqrt[q^n]{\theta}\vert_{\FF_q(X)_{\perf}^{\ab}}$ belongs to the image of the reciprocity law homomorphism $\rho$,
and hence $\theta$ has the first property.
\qed  $\;$\\

The following is our main result.
\begin{Theorem}\label{Theorem:MainResult}
Fix a nonzero effective divisor $D$ of $X$.
Also fix nonzero 
$$\alpha,\beta\in H^0(\OO_X(D)/\OO_X),$$ and a divisor $A_0$ of $X$ supported away from $D$ such that $$\deg A_0=g-2.$$ Then there exists a unique element $\varphi$ of the fraction field of $\Dbold$ such that for all $\theta\in \AUT{X}$, the following statements hold.
Firstly,
\begin{equation}\label{equation:BigDealInterpolation}
((1\otimes \theta)\varphi)\vert_\Delta=
\Hyp_D(\alpha,\beta,A_0+r_D(a))^{\min(\norm{a},1)}
\end{equation}
if $\theta$ is critical of exponent $a$ and the right side is defined.
Secondly,
\begin{equation}\label{equation:BigDealInterpolationBis}
((1\otimes \theta)\varphi)\vert_\Delta\neq 0,\infty
\end{equation}
if $\theta$ is not critical.
\end{Theorem}

Some amplifying remarks are in order.
\begin{enumerate}
\item Formula (\ref{equation:BigDealInterpolation}) is the  interpolation formula mentioned in the introduction. \item Formula (\ref{equation:BigDealInterpolationBis}) forces $\varphi$ to be a Coleman unit. See Prop.~\ref{Proposition:ColemanUnit} and its proof for a detailed explanation of this point.
\item Lemma~\ref{Lemma:UniquenessPrinciple}
already proves the uniqueness asserted in the theorem.
\item Lemma~\ref{Lemma:UniquenessPrinciple}, the theorem and relation (\ref{equation:MooreHyperTer})
among hypergeometric ratios in the high degree case force (\ref{equation:MooreHyperTer})  to hold  in the low degree case.
\item  Lemma~\ref{Lemma:Slick} simplifies the task of recognizing when $\theta$ is critical.
\item The theorem says nothing about $\varphi$ in the case that
$\theta$ is critical of exponent $a$ such that the right side of (\ref{equation:BigDealInterpolation}) is undefined---but the gap is filled by the author's  conjecture \cite[Conj.\ 9.5]{AndersonStirling}.
\end{enumerate}
In \S\ref{section:Discussion} we provide further amplification of the theorem,
in particular indicating the position of the theorem with respect to the author's conjecture.

The proof of the theorem commences in \S\ref{section:Tools}
and takes up the rest of the paper.
In \S\ref{section:Tools} we collect tools for the proof and in particular we put what we need of geometric class field theory into a form compatible with the Thakur-style approach to shtukas.
In \S\ref{section:CatalanDrinfeldSymbol} 
we study rank one shtukas and
we set up the Catalan-Drinfeld symbol formalism. A version  of Drinfeld's ``$\chi=0\Rightarrow h^0=h^1=0$'' \linebreak lemma (Lemma~\ref{Lemma:ChiZero} below) plays the key role. The Catalan-Drinfeld symbol formalism
is of intrinsic interest and no doubt further study of it will lead to
refinements of our conjecture.
In \S\ref{section:EndGame} we finish the proof
of Theorem~\ref{Theorem:MainResult}
by evaluating the Catalan-Drinfeld symbol in apt ways.

\section{Discussion}
\label{section:Discussion}

 We calculate hypergeometric ratios and verify  Theorem~\ref{Theorem:MainResult} ``by hand''
 in a simple special case.
We review the notion of Coleman unit and explain
why the functions produced by the theorem are
Coleman units.
We discuss the theorem in relation to  Coleman's paper \cite{Coleman}
and the author's conjecture \cite[Conj.\ 9.5]{AndersonStirling}.

\subsection{Sample calculation of hypergeometric ratios}\label{subsection:SimpleExample}

We assume under this heading that
$$X/\FF_q=\PP^1_t/\FF_q,\;\;\;\FF_q(X)=\FF_q(t).$$
 For each  $c\in \FF_q\cup\{\infty\}=\PP^1_t(\FF_q)$,
let $[c]$ be the corresponding
closed point of $\PP^1_t$. 
Let 
$$\alpha_\infty,\alpha_1,\alpha_0\in H^0(\OO_{\PP^1_t}([\infty]+[1]+[0])/\OO_{\PP^1_t})$$
be the $\FF_q$-basis consisting of elements represented by 
$$t,\frac{1}{1-t},\frac{t-1}{t}
\left(=\frac{1}{1-\frac{1}{1-t}}\right)\in H^0(\OO_{\PP^1_t}([\infty]+[1]+[0])),$$
respectively.  
We claim that
\begin{equation}\label{equation:BasicThreePoint}
\begin{array}{rcl}
\Hyp_{[\infty]+[0]}(\alpha_\infty,\alpha_0,(N-2)[1])&=&
t^{\epsilon_N\frac{q^{|N|}-1}{q-1}},\\
\Hyp_{[1]+[\infty]}(\alpha_1,\alpha_\infty,(N-2)[0])&=&
\left(\frac{1}{1-t}\right)^{\epsilon_N\frac{q^{|N|}-1}{q-1}},\\
\Hyp_{[0]+[1]}(\alpha_0,\alpha_1,(N-2)[\infty])&=&
\left(\frac{t-1}{t}\right)^{\epsilon_N\frac{q^{|N|}-1}{q-1}}
\end{array}
\end{equation}
for all nonzero integers 
$N$, where $\epsilon_N\in \{\pm 1\}$
is the sign of $N$. By symmetry (the map $t\mapsto 1/(1-t)$ is an automorphism of $\PP^1_t/\FF_q$ of order $3$) we have only to prove the first formula. Call the left side of the first formula $\Hyp(N)$ to abbreviate.
Note that the case $N>0$ (resp., $N<0$)  corresponds to the high
(resp., low) degree case of the definition of the hypergeometric ratio.
Assume at first that $N>0$. Take liftings $\tilde{\alpha}_\infty=t-1$ and $\tilde{\alpha}_0=\frac{(t-1)}{t}$. Then we have
$$\Hyp(N)=\prod_{e\in H^0(\OO_{\PP^1_t}((N-2)[1]))}
\frac{t-1+e}{\frac{(t-1)}{t}+e}$$
$$=\prod_{e\in H^0(\OO_{\PP^1_t}(-[0]+(N-1)[\infty]))}
\frac{t-1+\frac{e}{t(t-1)^{N-2}}}
{\frac{(t-1)}{t}+\frac{e}{t(t-1)^{N-2}}}$$
$$=\frac{\Moore(t^{N},t^{N-1},\dots,t)}{\Moore((-1)^{N-1},t^{N-1},\dots,t)}=\frac{\Moore(t^{N},t^{N-1},\dots,t)}{\Moore(t^{N-1},\dots,1)}=t^{\frac{q^{N}-1}{q-1}}.$$
We turn to the remaining case $N<0$. Put
$\nu=|N|$. We have
$$\Hyp(N)=\prod_{\begin{subarray}{c}\omega\in H^0(\Omega_{{\PP^1_t}/\FF_q}((\nu+2)[1]))\\
\RES_{[\infty]+[0]}(\frac{t-1}{t}\omega)=1
\end{subarray}} \omega \bigg/ \prod_{\begin{subarray}{c}\omega\in H^0(\Omega_{{\PP^1_t}/\FF_q}((\nu+2)[1]))\\
\RES_{[\infty]+[0]}(t\omega)=1
\end{subarray}} \omega
$$
$$=\prod_{\begin{subarray}{c}e\in H^0(\OO_{\PP^1_t}(\nu [\infty]))\\
\Res_{[0]}(t^{-1}e(t-1)^{-\nu-2}dt)=-1
\end{subarray}} e \bigg/ \prod_{\begin{subarray}{c}e\in H^0(\OO_{\PP^1_t}(\nu[\infty]))\\
\Res_{[\infty]}(te(t-1)^{-\nu-2}dt)=1
\end{subarray}} e$$
$$=\frac{\Moore((-1)^{\nu+1},t,\dots,t^\nu)}
{\Moore(t,\dots,t^\nu)}\bigg /
\frac{\Moore(-t^\nu,1,\dots,t^{\nu-1})}{\Moore(1,\dots,t^{\nu-1})}$$
$$=\frac{\Moore(1,\dots,t^{\nu-1})}{\Moore(t,\dots,t^\nu)}=
t^{-\frac{q^\nu-1}{q-1}}.$$
The claim is proved.

\subsection{Sample instance of theorem}\label{subsection:SimpleExampleBis}
Continuing in the setting of \S\ref{subsection:SimpleExample} above,
we verify Theorem~\ref{Theorem:MainResult} in the case 
\begin{equation}\label{equation:TheoremCase}
X=\PP_t^1,\;\;\;\FF_q(X)=\FF_q(t),\;\;\;(D,A_0,\alpha,\beta)=
([\infty]+[0],-2[1],\alpha_\infty,\alpha_0).
\end{equation}
Fix 
$$\tau\in \FF_q(X)^{\ab}\subset\overline{\FF_q(X)}$$ such that 
$$\tau^{q-1}=t.$$ 
Note that $[\infty]+[0]$ is a conductor for the abelian extension
$\FF_q(t,\tau)/\FF_q(t)$. 
Put $$\varphi=\tau^{-1}\otimes \tau\in 
\Dbold^\times.$$
We will verify that $\varphi$ has the  properties
(\ref{equation:BigDealInterpolation}) and (\ref{equation:BigDealInterpolationBis}) required by Theorem~\ref{Theorem:MainResult}.
Because $\varphi$ is a unit of $\Dbold$, condition (\ref{equation:BigDealInterpolationBis}) of the theorem is trivially satisfied by $\varphi$. Only condition (\ref{equation:BigDealInterpolation}) requires proof.
In more detail, what we need to prove is that 
\begin{equation}\label{equation:NeedToProveEx}
((1\otimes \theta)\varphi)\vert_\Delta=
\Hyp_{[\infty]+[0]}(\alpha_\infty,\alpha_0,-2[1]+r_{[0]+[\infty]}(a))^{\min(\norm{a},1)}
\end{equation}
for every $a\in \adeles_X^\times$ and $\theta\in \AUT{X}$ such that 
$$\rho^*(a)=\theta\vert_{\FF_q(X)_{\perf}^{\ab}},\;\;\;\;\norm{a}\neq 1.$$
Fix such $a$ and $\theta$ now, and also
fix $c\in \FF_q^\times$ and an integer $N\neq 0$ such that 
$$r_{[\infty]+[0]}(a)=
\left(\begin{array}{l}
\mbox{the generalized divisor class}\\
\mbox{of $[c]+(N-1)[1]$ of}\\
\mbox{conductor $[\infty]+[0]$}
\end{array}\right),\;\;\mbox{and hence}\;\;\norm{a}=q^N.$$ 
Since the image of $r_{[\infty]+[0]}$ is a copy of
$\ZZ\times\FF_q^\times$ we can indeed find $N$ and $c$
with these properties.
Notice now that $\rho(a)$ restricted to $\FF_q(t,\tau)$ is an arithmetic Frobenius element in $\Gal(\FF_q(t,\tau)/\FF_q(t))$ above $[c]$
(see the remark of \S\ref{subsubsection:Cancellation}) and hence
$$\rho(a)\tau=c\tau,\;\;\;\rho^*(a)\tau=c^{-1}\tau^{q^N}.$$
It follows that
\begin{equation}\label{equation:SimpleExampleBis}((1\otimes\theta)\varphi)\vert_\Delta=c^{-1}\tau^{q^N-1}=
\left\{\begin{array}{rl}
c^{-1}t^{\frac{q^N-1}{q-1}}&\mbox{if $N>0$,}\\
\sqrt[q^{|N|}]{c^{-1}t^{-\frac{q^{|N|}-1}{q-1}}}&\mbox{if $N<0$.}
\end{array}\right.
\end{equation}
Now by combining (\ref{equation:HypInv}) and (\ref{equation:MooreHyperBis})
with the first of the suite of formulas (\ref{equation:BasicThreePoint}), we have
\begin{equation}\label{equation:SimpleExample}
\Hyp_{[\infty]+[0]}(\alpha_\infty,\alpha_0,
[c]+(N-3)[1])=
c^{-1}t^{\epsilon_N\frac{q^{|N|}-1}{q-1}}.
\end{equation}
Finally, compare (\ref{equation:SimpleExampleBis}) and (\ref{equation:SimpleExample}) in order to see
that (\ref{equation:NeedToProveEx}) and hence 
(\ref{equation:BigDealInterpolation}) hold for $\varphi$.
The verification of the theorem in the special case (\ref{equation:TheoremCase}) is complete.

\subsection{The notion of Coleman unit}\label{subsection:ColemanUnits}
We review a  notion introduced in the author's paper \cite{AndersonStirling}
and inspired by Coleman's paper \cite{Coleman}. 
Definitions recalled under this heading are local to \S\ref{section:Discussion}
and will not be used  from \S\ref{section:Tools} onward.
\subsubsection{The ring  $\Kbold$}
Consider the subrings
$$
\widetilde{\Kbold}=\FF_q(X)^{\ab}_{\perf}\otimes_{\overline{\FF}_q}
\FF_q(X)^{\ab}_{\perf},\;\;\;\;\Kbold=\widetilde{\Kbold}^{\{\sigma\otimes \sigma\vert
\sigma\in \Gal(\FF_q(X)^{\ab}_{\perf}/\FF_q(X)_{\perf})\}}$$
of $\Dbold$.
By an evident modification of the proof of Lemma~\ref{Lemma:DGeography}, one verifies that for every maximal ideal $\Mbold\subset \Kbold$,
the corresponding local ring $\Kbold_\Mbold$ is a nondiscrete valuation ring
of rank $1$. Note that since $\Dbold/\Kbold$ is an integral extension
of domains, every maximal ideal of $\Kbold$ lies below some maximal
ideal of $\Dbold$. Note also that the image of diagonal evaluation  restricted to $\Kbold$ is $\FF_q(X)_{\perf}$.

\subsubsection{The twisting action}
We define the {\em twisting action}
$\varphi\mapsto \varphi^{(a)}$ of $\adeles_X^\times$ on $\widetilde{\Kbold}$ by the rule
$$(x\otimes y)^{(a)}=(\rho(a)x)\otimes y^{\norm{a}}.$$
Note that the twisting action stabilizes $\Kbold$.
We remark that for every $a\in \adeles^\times$ the automorphisms
$$x\otimes y\mapsto (x\otimes y)^{(a)},\;\;\;\;x\otimes y\mapsto x\otimes \rho^*(a)y$$
of $\widetilde{\Kbold}$ agree on $\Kbold$. We extend the twisting action to the
fraction field of $\Kbold$ in the unique possible way.

\subsubsection{Definition of  Coleman unit}
According to the definition \cite[\S 9.4]{AndersonStirling},
a {\em Coleman unit} $\varphi$ is a nonzero element
of the fraction field of $\Kbold$ such that for every maximal
ideal $\Mbold\subset \Kbold$, if $\varphi$ fails to be a unit
of the local ring $\Kbold_\Mbold$,
then $\Mbold$ is of the form
$$\Mbold=\ker\left((\varphi\mapsto \varphi^{(a)}\vert_\Delta):\Kbold\rightarrow \FF_q(X)_{\perf}\right)$$
for some $a\in \adeles_X^\times$.\\

\begin{Proposition}\label{Proposition:ColemanUnit}
Fix $D$, $\alpha$, $\beta$, 
and $A_0$ as in Theorem~\ref{Theorem:MainResult}, 
and let $\varphi$ satisfy the conclusion of the theorem. Then $\varphi$
is a Coleman unit.
\end{Proposition}
\proof By Lemma~\ref{Lemma:UniquenessPrinciple}
and property (\ref{equation:BigDealInterpolation}) we have
$$(\sigma\otimes \sigma)\varphi=\varphi$$ for all $\sigma\in \Gal(\overline{\FF_q(X)}/\FF_q(X)_{\perf})$ and 
$$(1\otimes \sigma)\varphi=\varphi$$
for all $\sigma \in \Gal(\overline{\FF_q(X)}/\FF_q(X)^{\ab}_{\perf})$.
Thus $\varphi$ belongs to the fraction field of $\Kbold$.
Now fix a maximal ideal $\Mbold\subset \Kbold$ 
such that
$$\Mbold\neq \ker((\psi\mapsto \psi^{(a)}\vert_\Delta):\Kbold\rightarrow
\FF_q(X)_{\perf})$$ for every $a\in \adeles^\times$.
By Lemma~\ref{Lemma:DGeography} and integrality of the ring extension $\Dbold/\Kbold$,  for some $\theta\in \AUT{X}$, we have 
 $$\Mbold=\Kbold\cap\ker((\psi\mapsto ((1\otimes \theta)\psi)\vert_\Delta):\Dbold\rightarrow\Dbold).$$
By hypothesis concerning $\Mbold$, the automorphism
$\theta$ cannot be critical,
and hence by (\ref{equation:BigDealInterpolationBis}) it follows
that $\varphi$ is a unit of the local ring $\Kbold_{\Mbold}$.
Therefore $\varphi$ is indeed a Coleman unit.
\qed

\subsection{Coleman's function on the product
of a Fermat curve with itself}

We return to the settings of \S\ref{subsection:SimpleExample}
and \S\ref{subsection:SimpleExampleBis}, assuming as before that 
$$X=\PP^1_t,\;\;\;\FF_q(X)=\FF_q(t).$$
Fix $\tau_0,\tau_1\in \FF_q(X)^{\ab}$ such that
$$\tau_0^{q-1}=t,\;\;\;\tau_1^{q-1}=1-t.$$
Then $(\tau_0,\tau_1)$ is a generic point of the Fermat
curve $x^{q-1}+y^{q-1}=1$ over $\FF_q$.
Put
$$\varphi=\tau_0\otimes \tau_0^{-1}+\tau_1\otimes \tau_1^{-1}-1\in\Dbold.$$
Then $\varphi$ is the function
on the product of two copies of the Fermat curve of degree $q-1$ over $\FF_q$ considered in Coleman's paper \cite{Coleman}.
Let $A_0$ be a divisor of $\PP^1_t$ supported away from $[\infty]+[1]+[0]$
such that
$$\deg A_0=-2,\;\;\;A_0\sim_{[\infty]+[0]}-2[1],\;\;A_0\sim_{[\infty]+[1]}-2[0].$$
By an evident modification of the calculation undertaken in \S\ref{subsection:SimpleExampleBis} which uses not only
the first but also the second of the three formulas (\ref{equation:BasicThreePoint}), we have
\begin{equation}\label{equation:NeedToProveExBis}
((1\otimes \theta)\varphi)\vert_\Delta=
\Hyp_{[\infty]+[1]+[0]}(\alpha_0+\alpha_1-\alpha_\infty,\alpha_\infty, A_0+r_{[0]+[\infty]}(a))^{\min(\norm{a},1)}
\end{equation}
for all $a\in \adeles^\times$ such that $\log_q\norm{a}\neq -1,0$
and $\theta\in \AUT{X}$ such that $\theta\vert_{\FF_q(X)^{\ab}_{\perf}}=\rho^*(a)$.
In other words, Coleman's function $\varphi$ makes (\ref{equation:BigDealInterpolation}) hold for suitable data $(D,\alpha,\beta,A_0)$,
and therefore by Theorem~\ref{Theorem:MainResult} 
and Proposition~\ref{Proposition:ColemanUnit} must be a Coleman unit.
But it is actually easy to verify that $\varphi$ is a Coleman unit
``by hand''.
Indeed, the divisor of $\varphi$ on the product of two copies of the Fermat curve can be worked out exactly, and that is exactly what Coleman did
in \cite{Coleman} in order to carry out his remarkable elementary
analysis of the Frobenius endomorphism of the Jacobian of the Fermat curve of degree $q-1$ over $\FF_q$. 
Theorem~\ref{Theorem:MainResult} says that
Coleman-like functions are not special or isolated---rather, they are ubiquitous.

\subsection{Position of the main result  with respect to the author's conjecture}
If one rewrites the right side of formula (\ref{equation:BigDealInterpolation})  in terms of the Catalan symbol defined in \cite{AndersonStirling} using formula (\ref{equation:CatalanSymbolSpecialization}), one sees that Theorem~\ref{Theorem:MainResult} 
confirms the author's conjecture \cite[Conj.\ 9.5]{AndersonStirling} ``asymptotically'', i.~e., for $\max(\norm{a},\norm{a}^{-1})$ large.
Further, our conjecture granted,
every Coleman unit it produces must be constructible by natural operations
from the Coleman units which Theorem~\ref{Theorem:MainResult} produces; this follows from remark 
\cite[\S9.6.4]{AndersonStirling}.  
The proof of our conjecture thus comes down
to a straightforward (if rather long and painstaking)
analysis of the examples produced by Theorem~\ref{Theorem:MainResult} using
the adelic theory of \cite{AndersonStirling}
and the local theory of \cite{AndersonLocalStirling}.
We will provide the details on another occasion.

\section{Toolkit}\label{section:Tools}
We review what we need of geometric class field theory.
We put the needed material in a form 
compatible with the statement of
Theorem~\ref{Theorem:MainResult} and
 the Thakur-style approach to shtukas. We proceed rapidly,
assuming that the reader is familiar with the standard reference \cite{Serre}.
Proposition~\ref{Proposition:ExplicitReciprocity} below summarizes the discussion.  Along the way we formulate a very special case of Bertini's theorem (Lemma~\ref{Lemma:Bertini}) needed as a technical tool.
Notation introduced here is in force for the rest of the paper.

\subsection{Expansion of the setting for Theorem~\ref{Theorem:MainResult}}

\subsubsection{The universal domain $W$}
We have previously chosen an algebraic closure
$\overline{\FF_q(X)}/\FF_q(X)$ and defined
$\overline{\FF}_q$ to be the algebraic closure of $\FF_q$ in $\overline{\FF_q(X)}$. We now fix an algebraically closed field $W$ containing $\overline{\FF}_q$ as a subfield.
Save for requiring $W$ to contain $\overline{\FF}_q$,
we choose $W$ independently of our
 previous choice of algebraic closure $\overline{\FF_q(X)}/\FF_q(X)$.  Elements of $W$ will sometimes be called {\em constants}. The field $W$ will play the role of a Weil-style universal domain. Later we will need $W$ to be large enough to permit construction of an embedding
 $\Dbold\rightarrow W$, 
but for the moment we make no assumption concerning the absolute transcendence degree of $W$, so that the conclusions we draw here will be valid for any algebraically closed field extending
$\overline{\FF}_q$. Given a ring $R$ between $\FF_q$ and $W$, let $R_{\perf}$ be the closure of $R$ in $W$ under the
extraction of $q^{th}$ roots. Given a field $K$ between $\FF_q$ and $W$, let $\overline{K}$ be the algebraic closure of $K$ in $W$,
let $K^{\sep}$ be the separable algebraic closure of $K$ in $\overline{K}$,
let $K^{\ab}$ be the abelian closure of $K$ in $K^{\sep}$,
let $K^{\ab}_{\perf}=(K^{\ab})_{\perf}$, and finally, let $$X_K/K =X\times_{\Spec(\FF_q)}\Spec(K)/\Spec(K)$$
be the base-change
of $X/\FF_q$. 

\subsubsection{Points and divisors  defined over the universal domain}
Abusing notation, we write $\overline{X}=X(W)=X_W$.
In other words, sometimes $\overline{X}$ denotes the set of \linebreak $W$-valued
points of $X$ and sometimes $\overline{X}$ denotes the $W$-scheme $X_W$.
In context this  usage should not cause confusion.
Correspondingly, we identify the
free abelian group generated by the set $\overline{X}$ with the divisor group of the curve $\overline{X}$. 
Given a divisor $D$ of $\overline{X}$, let $\supp D\subset \overline{X}$
be the support of $D$.
Given a divisor $D$ of $X$, let the divisor of $\overline{X}$ obtained by base-change
be again denoted by $D$.
Given $\xi\in \overline{X}$, let $\FF_q(\xi)$ be the subfield of $W$
generated over $\FF_q$ by the coordinates of $\xi$. We say that $\xi$ is {\em generic}
if $\FF_q(\xi)$ is an isomorphic copy of $\FF_q(X)$. We say that two generic points $\xi,\eta\in \overline{X}$ are {\em independent}
if the fields $\FF_q(\xi)$ and $\FF_q(\eta)$ are linearly disjoint over $\FF_q$,
in which case $\overline{\FF_q(\xi)}$ and
$\overline{\FF_q(\eta)}$ are linearly disjoint over $\overline{\FF}_q$.

\subsubsection{Generalized divisor classes defined over the universal domain}\label{subsubsection:Gen2} 
We adapt the definitions given in \S\ref{subsubsection:Gen1}
for $X/\FF_q$ to $\overline{X}/W$ in the obvious way.
\pagebreak
\begin{Lemma}\label{Lemma:Bertini}
Let $D$ be an effective divisor of $X$.
Let $E$ be a divisor of $\overline{X}$
supported away from $D$. Let $S$ be a finite subset of $\overline{X}\setminus \supp D$. 
Then there exist
\begin{itemize}
\item a divisor $\tilde{E}$ of $\overline{X}$, and
\item a divisor $A$ of $X$
\end{itemize}
such that:
\begin{itemize}
\item $\tilde{E}$ and $A$ are effective and supported away from $D$;
\item $\tilde{E}\sim_D A+E$;
\item $S\cap\supp \tilde{E}=\emptyset=S\cap \supp A$; and
\item $\tilde{E}$ has multiplicity $\leq 1$ everywhere on $\overline{X}$.
\end{itemize}
\end{Lemma}
\proof 
Let $\overline{X}_D$ be the singular model of
$\overline{X}$ constructed according to the procedure of \cite[Chap.\ IV, \S4]{Serre}.  
Roughly speaking, $\overline{X}_D$ is obtained by crushing $D$ to a single closed point $\infty_D$. 
Choose an effective divisor $A$ of $X$ of positive degree supported away from $D$ and $S$.
Consider the space
$$V=\{f\in H^0(\OO_{\overline{X}}(A+E))\mid f\vert_D\;\mbox{is constant}\},$$
and choose a $W$-basis $v_0,\dots,v_n\in V$.
After replacing $A$ by a sufficiently high multiple of itself,
we may assume that the map
$$(v_0:\cdots:v_n):\overline{X}\rightarrow \PP^n_{/W}$$
is a projective embedding of $\overline{X}_D$. For simplicity, let us identify
$\overline{X}_D$ with its image under this projective embedding,
and in turn identify $\overline{X}\setminus\supp D$
with $\overline{X}_D\setminus\{\infty_D\}$.
 Any hyperplane section $H\cap \overline{X}_D$ 
 to which $\infty_D$ does not belong
can then be construed as a member of
the generalized divisor class of $A+E$ of conductor $D$.
But any sufficiently general hyperplane $H$
does not intersect $S\cup\{\infty_D\}$ and by Bertini's
theorem cuts $\overline{X}_D$ transversely.
Take $\tilde{E}=H\cap \overline{X}_D$
for a general hyperplane $H$.
\qed

\subsection{Further expansion of the setting}
We prepare to state a version of explicit reciprocity.
\subsubsection{Twisting}
Given $\xi\in \overline{X}$ and an integer $n$ (possibly negative), let $\xi^{(n)}$ be the result of applying
the $(q^n)^{th}$ power automorphism of $W$ to $\xi$,  and let the map $(\xi\mapsto \xi^{(n)}):\overline{X}\rightarrow \overline{X}$ thus defined be extended additively to the group of divisors of $\overline{X}$.
Let $f\mapsto f^{(n)}$ be the unique automorphism of the function field of $\overline{X}$
which restricts on $W$ to the $(q^n)^{th}$ power automorphism of $W$ and which restricts
on the function field of $X$ to the identity automorphism.  
We call the operations $D\mapsto D^{(n)}$ on divisors and
$f\mapsto f^{(n)}$ on functions {\em $n$-fold twisting}.
Twisting commutes with formation 
of principal divisors, i.~e., $(f^{(n)})=(f)^{(n)}$.
 A meromorphic function $f$ on $\overline{X}$
satisfies $f^{(1)}=f$ if and only if $f$ descends to a meromorphic function on $X$.
Similarly, a divisor $D$ of $\overline{X}$ satisfies $D^{(1)}=D$
if and only if $D$ descends to a divisor of $X$.
For each effective divisor $D$ of $X$, $n$-fold
twisting preserves the group of divisors principal to the conductor $D$.

\subsubsection{Conjugation of divisors}
Given algebraically closed subfields $L_1,L_2\subset W$,
a divisor $E$ of $\overline{X}$ with $\supp E\subset X(L_1)$,
and an $\FF_q$-linear isomorphism \linebreak $\theta:L_1\iso L_2$,
let $\theta(E)$  be the result of applying the unique additive extension of the map  $(\zeta\mapsto \theta(\zeta)):X(L_1)\iso X(L_2)$ to $E$.   The operation $\theta$ fixes every divisor of $X$. If now further we are given an effective divisor $D$ of $X$, and we suppose
that $E$ is supported away from $D$ and satisfies $E\sim_D0$,
then necessarily $\theta(E)\sim_D0$.

\subsubsection{Key exact sequences}
Let $D$ be an effective divisor of $X$.
Let $J_D/\FF_q$ be the generalized Jacobian
of $X/\FF_q$ of conductor $D$, as defined by Rosenlicht. Put $\overline{J}_D=J_D(W)$.
Then $J_D(\FF_q)$  (resp., $\overline{J}_D$) is canonically equal to the group of generalized divisor classes of $X$ (resp., $\overline{X}$) of conductor $D$ and degree $0$. 
Crucially:
\begin{equation}\label{equation:Crucial}
\begin{array}{l}
\mbox{The operation
$E\mapsto E^{(1)}$ on
divisors of $\overline{X}$ supported away}\\
\mbox{from $D$ induces a map $\overline{J}_D\rightarrow \overline{J}_D$
equal to that induced by}\\
\mbox{the $q^{th}$ power Frobenius endomorphism $\Frob_q:J_D\rightarrow J_D$.}
\end{array}
\end{equation}
 We have an exact sequence
\begin{equation}\label{equation:LangTorsor}
0\rightarrow J_D(\FF_q)\subset
\overline{J}_D\xrightarrow{x\mapsto (1-\Frob_q)x}\overline{J}_D\rightarrow 0
\end{equation}
compatible with conjugation of divisors
at our disposal, due to Lang.
The preceding exact sequence is invariably applied
below in conjunction
with the exact sequence
\begin{equation}\label{equation:LangWeil}
0\rightarrow J_D(\FF_q)
\rightarrow
\left(\begin{array}{l}
\mbox{group of generalized divisor}\\
\mbox{classes of $X$ of conductor $D$}
\end{array}\right)
\xrightarrow{E\mapsto \deg E}\ZZ\rightarrow 0
\end{equation}
the existence of which is well-known
(for example, see  \cite[Cor.\ 4, Chap.\ VII, \S5]{Weil}).

\subsubsection{Slightly modified versions of $\rho$ and $\rho^*$}
Suppose we are given a generic point $\xi\in \overline{X}$. 
Let 
$$\rho_\xi:\adeles_X^\times\rightarrow \Gal(\FF_q(\xi)^{\ab}_{\perf}/\FF_q(\xi)_{\perf})$$ 
be the result of composing 
the reciprocity law homomorphism $\rho$ with the 
isomorphism 
$$\Gal(\FF_q(X)^{\ab}_{\perf}/\FF_q(X)_{\perf})\iso
\Gal(\FF_q(\xi)^{\ab}_{\perf}/\FF_q(\xi)_{\perf})$$
induced by the evaluation isomorphism 
$$(f\mapsto f\vert_\xi):\FF_q(X)\iso \FF_q(\xi).$$
Let
$$\rho_\xi^*:\adeles_X^\times\rightarrow\Aut(\FF_q(\xi)^{\ab}_{\perf}/\FF_q(\xi)_{\perf})$$
be defined by the rule
$$\rho^*_\xi(a)x=(\rho_\xi(a)^{-1}x)^{\norm{a}}$$
for all $x\in \FF_q(X)^{\ab}_{\perf}$.

\begin{Proposition}\label{Proposition:ExplicitReciprocity}
Fix a nonzero effective divisor $D$ of $X$. Fix a generic point $\xi\in \overline{X}$.
Fix  a divisor $I$ of $X$ supported away from $D$ of degree $1$. Fix a divisor $E$ of $\overline{X}$ supported away from $D$ such that
$$\deg E=0,\;\;\;E-E^{(1)}\sim_D \xi-I.$$
Fix $a\in \adeles_X^\times$, $N\in \ZZ$ and
$\mu\in \Aut(W/\overline{\FF}_q)$ such that
$$\norm{a}=q^N,\;\;\;\mu\vert_{\FF_q(\xi)^{\ab}_{\perf}}=\rho^*_\xi(a).$$
Then we have
\begin{equation}\label{equation:Explicit}
\mu(E)\sim_D E^{(1)}+ r_D(a)-I+
\left\{\begin{array}{cl}-\sum_{k=1}^{N-1}\xi^{(k)}&\mbox{if $N>1$,}\\
0&\mbox{if $N=1$,}\\
\sum_{k=0}^{|N|}
\xi^{(-k)}&\mbox{if $N\leq 0$.}
\end{array}\right.
\end{equation}
\end{Proposition}
\noindent Once we have proved the proposition we are free of any further necessity to discuss generalized Jacobians. Knowledge of the facts (\ref{equation:Crucial}--\ref{equation:Explicit}) will suffice.
We will be able to do all our work by manipulating divisors and functions on $\overline{X}$, just as in Thakur's paper \cite{Thakur}. For convenient application to the proof of Theorem~\ref{Theorem:MainResult} we have emphasized 
$\rho^*$ rather than $\rho$.
\proof We may assume without loss of generality
that $W=\overline{\FF_q(\xi)}$.
There exists a unique morphism
$$\Abel_{D,I}:X\setminus D\rightarrow J_D$$
of $\FF_q$-schemes 
such that
$$\Abel_{D,I}(\bar{x})=
\left(\begin{array}{l}
\mbox{generalized divisor class}\\
\mbox{of $\bar{x}-I$ of conductor $D$}
\end{array}\right)$$
for all points $\bar{x}\in \overline{X}\setminus \supp D$.
Put
$$\tau=\left(\begin{array}{l}
\mbox{generalized divisor class}\\
\mbox{of $E$ of conductor $D$}
\end{array}\right)\in \overline{J}_D.$$
Then $\tau$ is a solution of the Lang torsor equation
$$
(1-\Frob_q)(\tau)=\Abel_{D,I}(\xi).
$$
Now according to Lang we know that 
$$(1-\Frob_q):J_D\rightarrow J_D$$ is finite \'{e}tale surjective,
and hence $\tau\in J_D(K)$ for some finite
subextension $K/\FF_q(\xi)$ of $\FF_q(\xi)^{\ab}/\FF_q(\xi)$.
Now let a place $v$ of $\FF_q(\xi)$ unramified
in $K/\FF_q(\xi)$ be given and let $x$ be the closed point of $X$ corresponding to $v$ under the
isomorphism $(f\mapsto f\vert_\xi):\FF_q(X)\iso\FF_q(\xi)$.
Suppose that $\sigma_v\in \Aut(W/\FF_q(\xi)_{\perf})$ restricts to an arithmetic Frobenius element in $\Gal(K/\FF_q(\xi))$ at $v$.
Then we have
\begin{equation}\label{equation:ReciprocityIntermediate}
(1-\sigma_v)\tau=\left(\begin{array}{l}
\mbox{generalized divisor class of }\\
\mbox{$x-(\deg x)I$ of conductor $D$}
\end{array}\right).
\end{equation}
Now let $\sigma\in \Aut(W/\FF_q(\xi)_{\perf})$ be defined
by the rule
$$(\sigma(x))^{q^N}=\mu(x),$$
in which case
$$\sigma\vert_{\FF_q(\xi)^{\ab}_{\perf}}=\rho_\xi(a)^{-1}.$$
Then from (\ref{equation:ReciprocityIntermediate}) and the remark of \S\ref{subsubsection:Cancellation} we deduce that
$$(1-\rho_\xi(a))\tau=(\sigma-1)\tau=\left(\begin{array}{l}
\mbox{generalized divisor class of }\\
\mbox{$r_D(a)-NI$ of conductor $D$}
\end{array}\right),$$
or equivalently,
\begin{equation}\label{equation:ReciprocityIntermediateBis}
\sigma(E)\sim_D E+r_D(a)-NI.
\end{equation}
One verifies easily that
\begin{equation}\label{equation:ReciprocityIntermediateTer}
E^{(N)}\sim_D E^{(1)}+(N-1)I+\left\{\begin{array}{cl}-\sum_{k=1}^{N-1}\xi^{(k)}&\mbox{if $N>1$,}\\
0&\mbox{if $N=1$,}\\
\sum_{k=0}^{|N|}
\xi^{(-k)}&\mbox{if $N\leq 0$.}
\end{array}\right.
\end{equation}
Finally, apply the $N$-fold twisting operation to both sides
of (\ref{equation:ReciprocityIntermediateBis}),
and then apply (\ref{equation:ReciprocityIntermediateTer}) to obtain (\ref{equation:Explicit}).
\qed

 \section{Invariants of rank one shtukas}\label{section:CatalanDrinfeldSymbol}

We take a relatively elementary point of view on rank one shtukas similar to that taken in Thakur's paper \cite{Thakur}.  By means of a variant (Lemma~\ref{Lemma:ChiZero}) of
 Drinfeld's marvelous ``$\chi=0\Rightarrow h^0=h^1=0$'' lemma
   we prove a result (Theorem~\ref{Theorem:ShtukaCohomology})
 giving us control of the cohomology of shtukas.
We  apply the result to justify the definition of the {\em Catalan-Drinfeld symbol}.
We show how to realize all hypergeometric ratios as values
of the Catalan-Drinfeld symbol (Props.~\ref{Proposition:Realization1} and \ref{Proposition:Realization2}).  We also write out a determinantal formula (Prop.~\ref{Proposition:Realization3}) for the value of the Catalan-Drinfeld symbol.

\subsection{Shtukas}
\subsubsection{Definition}
We call a quadruple 
$$(D,\xi,\eta,E)$$
a {\em shtuka} under the following conditions:
\begin{itemize}
\item $D$ is a nonzero effective divisor of $X$.
\item $\xi,\eta\in \overline{X}\setminus \supp D$.
\item $E$ is a divisor of $\overline{X}$ supported away from $D$.
\item $E-E^{(1)}\sim_D -\xi^{(1)}+\eta$.
\item $\deg E=g-1$.
\end{itemize}
We call $D$, $\xi$, $\xi^{(1)}$, $\eta$ and $E$ the {\em conductor}, {\em basepoint}, {\em pole},
{\em zero},  and {\em divisor} of the shtuka $(D,\xi,\eta,E)$,
respectively. The notion of the basepoint of a shtuka has not previously been emphasized and here will be crucial.

\begin{Lemma}\label{Lemma:ShtukaEU}
(i) For all nonzero effective divisors $D$ of $X$ and points $\xi,\eta\in \overline{X}$
such that $\xi,\eta\not\in \supp D$, there exists a divisor $E$ of $\overline{X}$ supported away from $D$ such that
$(D,\xi,\eta,E)$ is a shtuka.  (ii) The generalized divisor class of $E$ of conductor $D$ is unique up to the addition of a generalized divisor class of $X$ of conductor $D$ and degree zero. 
\end{Lemma}
\proof Exact sequences (\ref{equation:LangTorsor})
and (\ref{equation:LangWeil}) prove this.
\qed

\subsubsection{Special functions attached to a shtuka}
By definition there is associated to the shtuka $(D,\xi,\eta,E)$ a unique meromorphic function $f_{D,\xi,\eta,E}$
on $\overline{X}$ such that
\begin{itemize}
\item $f_{D,\xi,\eta,E}\vert_D\equiv 1$, and 
\item $(f_{D,\xi,\eta,E})=E^{(1)}-E-\xi^{(1)}+\eta$. 
\end{itemize}
Were we to follow the terminology of \cite{Thakur} more closely, we would
actually call the function $f_{D,\xi,\eta,E}$ a shtuka, but we prefer not to do so. Note that if $\{\xi^{(1)},\eta\}\cap E=\emptyset$,
and $\xi^{(1)}\neq \eta$, then the pole and zero of the shtuka are in fact a pole and a zero of $f_{D,\xi,\eta,E}$, respectively.
Let 
$$\Psi(D,\xi,\eta,E)=
\left\{\psi\in H^0(\overline{X},\OO_{\overline{X}}(E+D))\left|
\begin{array}{l}
\mbox{$\psi-\psi^{(1)}$ is regular in some}\\
\mbox{neighborhood of $D$.}\end{array}\right.\right\}.$$
 Another way to describe $\Psi(D,\xi,\eta,E)$ is as
the subset of $H^0(\OO_{\overline{X}}(E+D))$ consisting
of liftings of  elements of $H^0(\OO_X(D)/\OO_X)$ via the exact sequence
$$0\rightarrow H^0(\OO_{\overline{X}}(E))
\rightarrow H^0(\OO_{\overline{X}}(E+D))\rightarrow
H^0(\OO_{\overline{X}}(D)/\OO_{\overline{X}}).$$
From the latter point of view it is clear that we have an exact sequence
\begin{equation}\label{equation:ExactPsi}
0\rightarrow
\Psi(D,\xi,\eta,E)\cap H^0(\OO_{\overline{X}}(E))
\rightarrow \Psi(D,\xi,\eta,E)\xrightarrow{\star_{D,\xi,\eta,E}} 
H^0(\OO_X(D)/\OO_X)
\end{equation}
at our disposal.

\subsubsection{Nondegeneracy}
We say that a shtuka $(D,\xi,\eta,E)$ is {\em nondegenerate} if the following condition holds:
\begin{itemize}
\item $\eta\not\in \{\xi^{(i)}\mid i\in [1-\deg(E+ D),+\infty)\cap(-\infty,g]\cap \ZZ\}$.
\end{itemize}
We remark that the set in question here is empty if and only if
$(g,\deg D)=(0,1)$ if and only if $\deg (E+D)=0$.

\begin{Theorem}\label{Theorem:ShtukaCohomology}
Let $(D,\xi,\eta,E)$ be a nondegenerate shtuka.
Then the following hold: 
\begin{enumerate}
\item  $h^i(\OO_{\overline{X}}(E))=0$ for $i=0,1$. 
\item $\Psi(D,\xi,\eta,E)\otimes_{\FF_q}W=H^0(\OO_{\overline{X}}(E+D))$.
\item $\Psi(D,\xi,\eta,E)\cap H^0(\OO_{\overline{X}}(E+D-\xi))=\{0\}$.
\end{enumerate}
\end{Theorem}
\noindent First we need a lemma.

\begin{Lemma}\label{Lemma:ChiZero}
 Let $(D,\xi,\eta,E)$ be any shtuka
and put $f=f_{D,\xi,\eta,E}$. Fix
\begin{itemize}
\item  a divisor $D_1$ of $X$,
\item a nonnegative integer $m$, and
\item a positive integer $N$
\end{itemize}
such that
\begin{itemize}
\item $m\leq \deg (E+D_1)$, and
\item $\eta\not\in \{\xi^{(i)}\vert i\in [1-\deg (E+D_1),N-m]\cap \ZZ\}$.
\end{itemize}
Let there be given
$$0\neq \psi\in H^0(\OO_{\overline{X}}(E+D_1-\sum_{i=0}^{m-1}\xi^{(-i)})),$$
and, for every integer $k\geq 0$, 
define $\psi_k$ by the rules
$$\psi_0=\psi,\;\;\;\psi_{k+1}=f\psi_k^{(1)}.$$
Then the functions $\psi_0,\dots,\psi_{N}$ are $W$-linearly independent.
\end{Lemma}
\noindent This is a refinement of \cite[Lemma 3.3.1]{AndersonAHarmonic} and a direct descendant of Drinfeld's ``$\chi=0\Rightarrow h^0=h^1=0$'' lemma.
For the latter see \cite{Drinfeld} or \cite[p.\ 146]{Mumford}.  
\proof After replacing $m$ by a larger integer if necessary, we may assume without loss of generality that
\begin{equation}\label{equation:RefinedHypothesis}
\psi\in H^0(\OO_{\overline{X}}(E+D_1-\sum_{i=0}^{m-1}\xi^{(-i)}))\setminus H^0(\OO_{\overline{X}}(E+D_1-\sum_{i=0}^{m}\xi^{(-i)})).
\end{equation}
Put 
$$E_1=E-\sum_{i=0}^{m-1}\xi^{(-i)},\;\;\;\xi_1=\xi^{(-m)},$$
noting that
\begin{equation}\label{equation:Nifty}
(f)=E_1^{(1)}-E_1-\xi_1^{(1)}+\eta.
\end{equation}
For $k\geq -1$ put
$$\Xi_k=\left\{\begin{array}{rl}
-\xi_1&\mbox{if $k=-1$,}\\
0&\mbox{if $k=0$,}\\
\sum_{i=1}^k \xi_1^{(i)}&\mbox{if $k>0$,}
\end{array}\right.$$
noting that
\begin{equation}\label{equation:NiftyBis}
\Xi_k=\Xi_{k-1}^{(1)}+\xi_1^{(1)}=\Xi_{k-1}+\xi_1^{(k)}\;\;\mbox{for $k\geq 0$}.
\end{equation}
We claim that
$$\psi_k\in H^0(\OO_{\overline{X}}(E_1+D_1+\Xi_k))\setminus H^0(\OO_{\overline{X}}(E_1+D_1+\Xi_{k-1}))\;\;\mbox{for $k=0,\dots,N$.}$$
 The case $k=0$ is our hypothesis (\ref{equation:RefinedHypothesis}). For $N\geq k>0$, we have
$$\psi_k=f\psi_{k-1}^{(1)}\in H^0(\OO_{\overline{X}}(E_1+D_1+\Xi_k-\eta))\setminus H^0(\OO_{\overline{X}}(E_1+D_1+\Xi_{k-1}-\eta))$$
by (\ref{equation:Nifty},\ref{equation:NiftyBis}) and induction on $k$, and we have
$$\eta\neq\xi^{(k-m)} =\xi_1^{(k)}$$ by hypothesis,  so the claim holds in general.
The claim granted, the lemma is proved. \qed
\proof[Proof of Theorem~\ref{Theorem:ShtukaCohomology}]
Put $f=f_{D,\xi,\eta,E}$.

(i) Supposing that statement (i) fails, we have $g>0$ and we can find
some 
$$0\neq \psi\in H^0(\OO_{\overline{X}}(E)).$$
The lemma in the case 
$$(D_1,m,N)=(0,0,g)$$
combined with our hypothesis of nondegeneracy yields 
$W$-linearly independent functions 
$$\psi=\psi_0,\dots,\psi_{g}\in H^0(\OO_{\overline{X}}(E+\sum_{i=1}^g \xi^{(i)})).$$
But
$$\deg (E+\sum_{i=1}^g \xi^{(i)})=2g-1>2g-2$$
and hence
$$h^0(\OO_{\overline{X}}(E+\sum_{i=1}^g \xi^{(i)}))=g.$$ 
This contradiction proves statement (i).

(ii) By statement (i),
the natural sequence
$$0=H^0(\OO_{\overline{X}}(E))
\rightarrow H^0(\OO_{\overline{X}}(E+D))\rightarrow
H^0(\OO_{\overline{X}}(D)/\OO_{\overline{X}})\rightarrow 0$$
is exact, hence the homomorphism $\star_{D,\xi,\eta,E}$ in exact sequence (\ref{equation:ExactPsi})  is an isomorphism, and hence statement (ii) holds.

(iii) Supposing now that statement (iii) fails, 
there exists 
$$0\neq \psi\in \Psi(D,\xi,\eta,E)
\cap H^0(\OO_{\overline{X}}(E+D-\xi)).$$
Since the case $(g,\deg D)=(0,1)$ is already ruled out, we have
$\deg (E+D)>0$.
The lemma 
in the case 
$$(D_1,m,N)=(D,1,1)$$ 
combined with our hypothesis of nondegeneracy produces $W$-linearly independent functions
$$\psi_0=\psi\in H^0(\OO_{\overline{X}}(E+D-\xi)),\;\;
\psi_1=f\psi^{(1)}\in H^0(\OO_{\overline{X}}(E+D-\eta))$$
from which, since $f\vert_D\equiv 1$, we get a nonzero function
$$\psi_1-\psi_0\in H^0(\OO_{\overline{X}}(E)).$$
But the latter space is $0$-dimensional by statement (i). This contradiction proves statement (iii).
\qed

\subsection{The Catalan-Drinfeld symbol}
Let $(D,\xi,\eta,E)$ be a nondegenerate shtuka. For each 
$$\alpha\in H^0(\OO_X(D)/\OO_X)$$
there exists by Theorem~\ref{Theorem:ShtukaCohomology} 
a unique lifting
$$\psi_\alpha\in \Psi(D,\xi,\eta,E)$$ with respect to the exact
sequence
\begin{equation}\label{equation:ExactitudeBis}
0=H^0(\OO_{\overline{X}}(E)
\rightarrow H^0(\OO_{\overline{X}}(E+D))
\rightarrow H^0(\OO_{\overline{X}}(D)/\OO_{\overline{X}})\rightarrow 0,
\end{equation}
and moreover for $\alpha\neq 0$, the order of vanishing of
the meromorphic function
$\psi_\alpha$ at the point $\xi$ is independent of $\alpha$.
For all nonzero $\alpha,\beta\in H^0(\OO_X(D)/\OO_X)$ we define
$$\left[\begin{array}{cccc}
D&\xi&\eta&E\\
&\alpha&\beta&
\end{array}\right]=(\psi_\alpha/\psi_\beta)(\xi)\in W^\times,$$
which depends only on the generalized divisor class of $E$ to the conductor $D$.  We call the rule $\left[\begin{array}{cccc}
\cdot&\cdot&\cdot&\cdot\\
&\cdot&\cdot
\end{array}\right]$ the {\em Catalan-Drinfeld symbol}.
\begin{Proposition}\label{Proposition:Realization1}
Let $\xi\in \overline{X}$ be a generic point.
Let $N>g$ be an integer. Let $(D,\xi,\xi^{(N)},E)$ be a shtuka.
Then the following hold:
\begin{enumerate}
\item $(D,\xi,\xi^{(N)},E)$ is nondegenerate. 
\item There exists a divisor $E_0$ of $X$ supported away from $D$ such that
$$E\sim_D E_0-(\xi^{(1)}+\cdots+\xi^{(N-1)}).$$
\item We have
$$\left[\begin{array}{cccc}
D&\xi&\xi^{(N)}&E\\
&\alpha&\beta&
\end{array}\right]=\Hyp_D(\alpha,\beta,E_0)\vert_\xi
$$
for all nonzero $\alpha,\beta\in H^0(\OO_X(D)/\OO_X)$.
\end{enumerate}
\end{Proposition}
\proof Statement (i) is immediate. Statement (ii) follows from the definitions via exact sequence (\ref{equation:LangTorsor}).
We have only to prove statement (iii).
Without loss of generality we may assume that
$$E=E_0-(\xi^{(1)}+\cdots+\xi^{(N-1)}).$$
 By hypothesis $\deg E_0>2g-2$, hence $\Hyp_D(\alpha,\beta,E_0)$ is defined and moreover
$$h^0(\OO_X(E_0))=N-1,\;\;\;h^1(\OO_X(E_0))=0.$$
Choose an $\FF_q$-basis $e_1,\dots, e_{N-1}\in H^0(\OO_X(E_0))$.
Via the natural exact sequence
$$0\rightarrow H^0(\OO_X(E_0))\rightarrow
H^0(\OO_X(E_0+D))\rightarrow H^0(\OO_X(D)/\OO_X)\rightarrow 0$$
choose 
a lifting $\tilde{\alpha}\in H^0(\OO_X(E_0+D))$ of $\alpha$.
Since $h^i(\OO_{\overline{X}}(E))=0$ for $i=0,1$ by Theorem~\ref{Theorem:ShtukaCohomology},
there exist unique constants $C_1,\dots,C_{N-1}\in W$ such that
$$\psi_\alpha=\tilde{\alpha}-\sum_{i=1}^{N-1}C_ie_i\in H^0(\OO_{\overline{X}}(E+D)).$$
Put $C_N=\psi_\alpha(\xi)$. The coefficients $C_1,\dots,C_N$ satisfy the matrix equation
$$
\left[\begin{array}{cccc}
e_{1}(\xi^{(N-1)})&\dots&e_{N-1}(\xi^{(N-1)})&0\\
\vdots&&\vdots&\vdots\\
e_{1}(\xi^{(1)})&\dots&e_{N-1}(\xi^{(1)})&0\\
e_{1}(\xi^{(0)})&\dots&e_{N-1}(\xi^{(0)})&1
\end{array}\right]
\left[\begin{array}{c}
C_1\\
\vdots\\
C_{N-1}\\
C_{N}
\end{array}\right]=\left[\begin{array}{c}
\tilde{\alpha}(\xi^{(N-1)})\\
\vdots\\
\tilde{\alpha}(\xi^{(1)})\\
\tilde{\alpha}(\xi^{(0)})\end{array}\right].
$$
By Cramer's Rule we have
$$\psi_\alpha(\xi)=C_N=
\left.(-1)^{N-1}\frac{\Moore(\tilde{\alpha},e_{1},\dots,e_{N-1})}{\Moore(e_{1},\dots,e_{N-1})^q}\right|_{\xi},$$ 
whence the claimed formula via the Moore determinant identity. \qed

\begin{Proposition}\label{Proposition:Realization2}
Let $\xi\in \overline{X}$ be a generic point.
Let $N>-2+g+\deg D$ be an integer.
 Let $(D,\xi^{(N)},\xi,E)$ be a shtuka.
Then the following hold:
\begin{enumerate} \item $(D,\xi^{(N)},\xi,E)$ is nondegenerate. 
\item
There exists a divisor $E_0$ of $X$ supported away from $D$
such that
$$E\sim_D E_0+\xi^{(0)}+\cdots+\xi^{(N)}.$$
\item We have
$$\left[\begin{array}{cccc}
D&\xi^{(N)}&\xi&E\\
&\alpha&\beta&
\end{array}\right]=\Hyp_D(\alpha,\beta,E_0)\vert_\xi
$$
for all nonzero $\alpha,\beta\in H^0(\OO_X(D)/\OO_X)$.
\end{enumerate}
\end{Proposition}
\proof  As in the proof of the preceding proposition, statements (i) and (ii) are easy to check. We 
have only to prove statement (iii). 

We pause to introduce some notation.
Given a meromorphic differential $\omega$ on $\overline{X}$, let $\RES_D \omega$ be the sum of the residues $\Res_{\bar{x}}\omega\in W$ extended over closed points
$\bar{x}$ of $\overline{X}$ in the support of $D$. If $D$ descends to a divisor
of $X$ and $\omega$ descends to a meromorphic differential on $X$, then $\RES_D\omega$ as defined
here coincides with $\RES_D\omega$ as previously defined in \S\ref{subsubsection:RES}.

We turn to the proof of statement (iii). 
We may assume without loss of generality that 
$$E= E_0+\xi^{(0)}+\cdots+\xi^{(N)}.$$
Fix a lifting $\tilde{\alpha}$
of $\alpha$ to a meromorphic function on $X$.
By hypothesis \linebreak $\deg E_0<-\deg D$, hence $\Hyp_D(\alpha,\beta,E_0)$ is defined, moreover
$$h^0(\Omega_X(-E_0))=N+1,\;\;\;h^1(\Omega_X(-E_0-D))=0,$$
and hence we can find an $\FF_q$-basis $\omega_0,\dots,\omega_N\in H^0(\Omega_X(-E_0))$
such that 
$$\RES_D \tilde{\alpha}\omega_k=\delta_{0k}\;\;\mbox{for $k=0,\dots,N$}.$$ 
Fix a nonzero meromorphic differential $\zeta$
on $X$ arbitrarily. By ``sum-of-residues-equals-zero''
 we have
 $$\sum_{i=0}^N(\omega_k/\zeta)^{q^i}\vert_{\xi}
 \Res_{\xi^{(i)}}\psi_\alpha\zeta=\sum_{i=0}^N \Res_{\xi^{(i)}}\psi_\alpha\omega_k =
 -\RES_D\psi_\alpha\omega_k=-\delta_{0k},
 $$
and hence, equivalently,
$$\left[\begin{array}{ccc}
(\omega_0/\zeta)^{q^N}\vert_\xi&\dots&
(\omega_N/\zeta)^{q^N}\vert_\xi\\
\vdots&&\vdots\\
(\omega_0/\zeta)\vert_\xi&\dots&
(\omega_N/\zeta)\vert_\xi
\end{array}\right]
\left[\begin{array}{c}
\Res_{\xi^{(N)}}\psi_\alpha\zeta\\
\vdots\\
\Res_{\xi^{(0)}}\psi _\alpha\zeta
\end{array}\right]
=
\left[\begin{array}{r}
-1\\
0\\
\vdots\\
0\end{array}\right].
$$
By Cramer's Rule we have
$$\Res_{\xi^{(N)}}\psi_\alpha\zeta=-
\left.\frac{\Moore(\omega_1/\zeta,\dots,\omega_N/\zeta)}
{\Moore(\omega_0/\zeta,\dots,\omega_N/\zeta)}\right|_\xi,
$$
whence the desired result now via the Moore determinant identity.
\qed
\subsection{A determinantal formula for the Catalan-Drinfeld symbol}
Fix a nondegenerate shtuka $(D,\xi,\eta,E)$
and nonzero $\alpha,\beta\in H^0(\OO_X(D)/\OO_X)$.
Suppose we can write 
$$E=E_1-E_2$$ where 
\begin{itemize}
\item $E_1$ and $E_2$ are supported away from $D$,
\item $E_2$ is effective and of multiplicity $\leq 1$ everywhere on $\overline{X}$,
\item The sets $\supp E_1$, $\supp E_2$ and $\{\xi\}$ are disjoint.
\end{itemize}
Put
$$n=\deg E_2,\;\;\;E_2=\sum_{i=1}^n \xi_i\;\;\;(\xi_i\in \overline{X}),\;\;\;
\xi_0=\xi.$$
We have at our disposal a natural exact sequence
\begin{equation}\label{equation:Exactitude}
0=H^0(\OO_{\overline{X}}(E_1-E_2))
\rightarrow H^0(\OO_{\overline{X}}(E_1))
\rightarrow H^0(\OO_{\overline{X}}/\OO_{\overline{X}}(-E_2))\rightarrow 0.
\end{equation} 
It follows in particular that
$$h^0(\OO_{\overline{X}}(E_1))=n,\;\;\;
h^1(\OO_{\overline{X}}(E_1))=0.$$
Choose any $W$-basis
$$f_1,\dots,f_n\in H^0(\OO_{\overline{X}}(E_1)).$$
We have at our disposal a natural exact sequence
$$
0\rightarrow H^0(\OO_{\overline{X}}(E_1))
\rightarrow H^0(\OO_{\overline{X}}(E_1+D))\rightarrow
H^0(\OO_{\overline{X}}(D)/\OO_{\overline{X}})\rightarrow 0.
$$
Choose any liftings
$$\tilde{\alpha},\tilde{\beta}\in H^0(\OO_{\overline{X}}(E_1+D))$$
of $\alpha$ and $\beta$, respectively. Put 
$$\mbox{$g_i=f_i$ for $i=1,\dots,n$, $f_0=\tilde{\alpha}$ and $g_0=\tilde{\beta}$.}$$
Note that for $i,j=0,\dots,n$, both $f_i$ and $g_i$ have no pole at $\xi_j$.
\begin{Proposition}\label{Proposition:Realization3}
Notation and hypotheses as above,
$$\det_{i,j=0}^n g_i(\xi_j)\cdot \left[\begin{array}{cccc}
D&\xi&\eta&E\\
&\alpha&\beta
\end{array}\right]=
\det_{i,j=0}^n f_i(\xi_j),$$
and moreover neither of the determinants vanish.
\end{Proposition}
\proof By exactness of (\ref{equation:Exactitude}) 
and distinctness of the points $\xi_1,\dots,\xi_n$, we have
$$\det_{i,j=1}^n f_i(\xi_j)=\det_{i,j=1}^n g_i(\xi_j)\neq 0.$$
Applying Cramer's Rule again,
as in the proof of Proposition~\ref{Proposition:Realization1},
we find that
 $$\psi_\alpha(\xi)=
 \det_{i,j=0}^nf_i(\xi_j)\bigg/\det_{i,j=1}^n f_i(\xi_j),$$
where $\psi_\alpha\in H^0(\OO_{\overline{X}}(E+D))$
is the unique lifting of $\alpha$ via exact sequence
(\ref{equation:ExactitudeBis}).
Moreover, since $\xi\not\in\supp(E_1-E_2+D)$,
we have $\psi_\alpha(\xi)\neq 0$
by Theorem~\ref{Theorem:ShtukaCohomology}.
Our conclusions for $\alpha$ hold for $\beta$ also.
The result follows.
\qed

\section{Proof of the main result}\label{section:EndGame}

\subsection{Reduction to a calculation of Catalan-Drinfeld symbols}
\subsubsection{Data for the theorem}
Let
$D$, $\alpha$, $\beta$ and $A_0$ be as specified in Theorem~\ref{Theorem:MainResult}. We also fix $\theta\in \AUT{X}$ arbitrarily, save for imposing without loss of generality
the following condition:  if $\theta$ is critical of exponent $a\in \adeles_X^\times$, then 
$$g-2+\log_q\norm{a}=\deg(A_0+r_D(a))\not\in [-\deg D,\infty)\cap(-\infty,2g-2].$$
The latter is precisely the condition under which $\Hyp_D(\alpha,\beta,A_0+r_D(a))$ is defined
should $\theta$ happen to be critical of exponent $a$.

\subsubsection{Embeddings}
We fix independent generic points $\xi,\eta\in \overline{X}$. 
(And so at this point we are imposing the further condition on $W$ that the latter
be of absolute transcendence degree at least $2$.)
Fix an $\overline{\FF}_q$-linear isomorphism
$$\lambda:\overline{\FF_q(\xi)}\iso\overline{\FF_q(X)}$$
such that
$$f=\lambda(f\vert_\xi)$$
for all $f\in \FF_q(X)$.
Let 
$$\iota:\overline{\FF_q(\xi)}\iso\overline{\FF_q(\eta)}$$
be an $\overline{\FF}_q$-linear isomorphism
such that
$$\iota(\xi)=\eta.$$
Let
$$\epsilon:(\mbox{compositum in $W$ of $\overline{\FF_q(\xi)}$
and $\overline{\FF_q(\eta)}$})
\iso(\mbox{fraction field of $\Dbold$})$$
be the unique isomorphism such that
$$\epsilon(x\iota(y))=
\lambda(x)\otimes \lambda(y)$$
for all $x,y\in \overline{\FF_q(\xi)}$.
Let 
$$\mu:\overline{\FF_q(\xi)}\iso\overline{\FF_q(\xi)}$$
be the unique $\overline{\FF}_q$-linear automorphism such that
$$\theta\lambda=\lambda \mu.$$
Then we have
\begin{equation}\label{equation:McGuffin}
((1\otimes \theta)(\epsilon(x\,\iota(y))))\vert_\Delta=\lambda(x\,\mu(y))
\end{equation}
for all $x,y\in \overline{\FF_q(\xi)}$. 

\subsubsection{The reduction}
We fix a divisor $I$ supported away from $D$
such that $\deg I=1$. We select a divisor $E$ of $\overline{X}$ supported away from $D$
such that
$$\supp E\subset X(\overline{\FF_q(\xi)}),\;\;\;
\deg E=0,\;\;\;E-E^{(1)}\sim_D\xi-I,$$
as is evidently possible by applying (\ref{equation:Crucial}) and (\ref{equation:LangTorsor}) with $W=\overline{\FF_q(\xi)}$.
One verifies  that
$$(D,\xi,\eta,A_0+I-E^{(1)}+\iota(E)),\;\;\;
(D,\xi,\mu(\xi),A_0+I-E^{(1)}+\mu(E))$$
are nondegenerate shtukas, immediately in the former case
since $\xi$ and $\eta$ are independent,
and via Lemma~\ref{Lemma:Slick} in the latter case.
Further, by Propositions~\ref{Proposition:ExplicitReciprocity},
\ref{Proposition:Realization1} and \ref{Proposition:Realization2} we have
$$\begin{array}{cl}
&\displaystyle\left[\begin{array}{cccc}
D&\xi&\mu(\xi)&A_0+I-E^{(1)}+\mu(E)\\
&\alpha&\beta
\end{array}\right]\\\\
=&\displaystyle\left.\begin{array}{ll}
\left(\Hyp_D(\alpha,\beta,A_0+r_D(a))\vert_{\xi}\right)^{\min(\norm{a},1)}&\mbox{if $\theta$ is critical of exponent $a$.}
\end{array}\right.
\end{array}
$$
It will therefore be enough to show that
\begin{equation}\label{equation:McGuffin1}
\varphi=\epsilon
\left[\begin{array}{cccc}
D&\xi&\eta&A_0+I-E^{(1)}+\iota(E)\\
&\alpha&\beta
\end{array}\right]\;\;\mbox{is defined,}
\end{equation}
\begin{equation}\label{equation:McGuffin0}
\varphi_0=\lambda\left[\begin{array}{cccc}
D&\xi&\mu(\xi)&A_0+I-E^{(1)}+\mu(E)\\
&\alpha&\beta
\end{array}\right]\;\;\mbox{is defined, and}
\end{equation}
\begin{equation}
\label{equation:McGuffin2}
((1\otimes \theta)\varphi)\vert_\Delta=\varphi_0\end{equation}
in order to
finish the proof of Theorem~\ref{Theorem:MainResult}.

\begin{Lemma}\label{Lemma:McGuffin}
There exist 
\begin{itemize}
\item divisors $E_1$ and $E_2$ of $\overline{X}$, and
\item a divisor $A_1$ of $X$ 
\end{itemize}
such that:
\begin{enumerate}
\item $A_1$, $E_1$ and $E_2$ are effective and supported away from $D$;
\item $\supp E_1\cup  \supp E_2\subset X(\overline{\FF_q(\xi)})$;
\item $E_2$ is of multiplicity $\leq 1$ everywhere on $\overline{X}$;
\item $A_0+I-E^{(1)}+\mu(E)\sim_D A_1-E_1-\mu(E_2)$;
\item $A_0+I-E^{(1)}+\iota(E)\sim_D A_1-E_1-\iota(E_2)$; 
\item the sets $\supp (A_1-E_1)$, $\supp \mu(E_2)$ and $\{\xi\}$
are disjoint; and
\item the sets $\supp (A_1-E_1)$, $\supp \iota(E_2)$ and $\{\xi\}$
are disjoint.
\end{enumerate}
\end{Lemma}
\proof
By Lemma~\ref{Lemma:Bertini} (applied with $W=\overline{\FF_q(\xi)}$) we can find
\begin{itemize}
\item a divisor $E_3$ of $\overline{X}$, and
\item a divisor $A_3$ of $X$
\end{itemize}
such that:
\begin{itemize}
\item $A_3$ and $E_3$ are supported away from $D$;
\item $E_3$ is effective;
\item $E_3\sim_DA_3-\mu(E)$;
\item $\supp E_3\subset X(\overline{\FF_q(\xi)})\setminus\{\xi\}$; and
\item $E_3$ has multiplicity $\leq 1$ everywhere on $\overline{X}$.
\end{itemize}
By Lemma~\ref{Lemma:Bertini} (again applied with $W=\overline{\FF_q(\xi)}$) we can find
\begin{itemize}
\item a divisor $E_1$ of $\overline{X}$, and
\item a divisor $A_1$ of $X$
\end{itemize}
such that:
\begin{itemize}
\item $A_1$ and $E_1$ are effective and supported away from $D$;
\item $E_1\sim_D A_1-(A_0+I-E^{(1)}+A_3)$; 
\item $\supp A_1\cap \supp E_3=\emptyset$; and
\item $\supp E_1\subset X(\overline{\FF_q(\xi)})\setminus(\{\xi\}\cup \supp E_3)$.
\end{itemize}
Then $A_1$, $E_1$ and $E_2=\mu^{-1}(E_3)$  have all
the desired properties. 
\qed

\subsection{Completion of the proof of Theorem~\ref{Theorem:MainResult}} Let $A_1$, $E_1$ and $E_2$ be as provided by Lemma~\ref{Lemma:McGuffin}.
We now apply Proposition~\ref{Proposition:Realization3}.
Put
$$K=\overline{\FF_q(\xi)},\;\;\;n=\deg E_2,\;\;\;
E_2=\sum_{i=1}^n\xi_i,\;\;\xi_0=\xi.$$
Choose a $K$-basis
$$f_1,\dots,f_n\in H^0(\OO_{X_K}(A_1-E_1))
\subset K\otimes_{\FF_q}H^0(\OO_X(A_1))$$
and liftings
$$\tilde{\alpha},\tilde{\beta}\in H^0(\OO_{X_K}(A_1-E_1+D))
\subset K\otimes_{\FF_q}H^0(\OO_X(A_1+D)).$$
Put
$$\mbox{$g_i=f_i$ for $i=1,\dots,n$, $f_0=\tilde{\alpha}$
and $g_0=\tilde{\beta}$.}$$
Then we have formulas
$$\det_{i,j=0}^n g_i(\iota(\xi_j))\cdot \left[\begin{array}{cccc}
D&\xi&\eta&A_0+I-E^{(1)}+\iota(E)\\
&\alpha&\beta
\end{array}\right]=
\det_{i,j=0}^n f_i(\iota(\xi_j)),
$$
$$\det_{i,j=0}^n g_i(\mu(\xi_j))\cdot\left[\begin{array}{cccc}
D&\xi&\mu(\xi)&A_0+I-E^{(1)}+\mu(E)\\
&\alpha&\beta
\end{array}\right]=
\det_{i,j=0}^n f_i(\mu(\xi_j))
$$
which verify (\ref{equation:McGuffin1},\ref{equation:McGuffin0}) and, in view of (\ref{equation:McGuffin}),
also prove (\ref{equation:McGuffin2}). The proof of Theorem~\ref{Theorem:MainResult} is complete.  \qed

\end{document}